\documentclass{aims}
\usepackage{amsmath}
  \usepackage{paralist}
  \usepackage{graphics} 
  \usepackage{epsfig} 
 \usepackage[colorlinks=true]{hyperref}
\hypersetup{urlcolor=blue, citecolor=red}

\def\currentvolume{X}
 \def\currentissue{X}
  \def\currentyear{200X}
   \def\currentmonth{XX}
    \def\ppages{X--XX}

\newcommand{\R}{\mathbb{R}}
\newcommand{\N}{\mathbb{N}}
\newcommand{\B}{\mathbb{B}}
\newcommand{\xb}{\bar{x}}

\newcommand{\ub}{\bar{u}}

\newcommand{\co}{\mbox{co }}
\renewcommand{\ae}{\mbox{a.e.} }

\newcommand{\abs}[1]{\left\vert#1\mathstrut\right\vert}
\newcommand{\disp}{\displaystyle}

\newcommand{\beq}{\begin{equation}}
\newcommand{\eeq}{\end{equation}}
\newcommand{\ba}{\begin{array}}
\newcommand{\ea}{\end{array}}
\newcommand{\beeq}{\begin{eqnarray}}
\newcommand{\eeeq}{\end{eqnarray}}

\renewcommand{\d}{\partial}
\newcommand{\norm}[1]{\left\Vert#1\mathstrut\right\Vert}
\newcommand{\Set}[1]{\left\lbrace#1\mathstrut\right\rbrace}
\newcommand{\lr}[3]{\left#1{\mathstrut#3}\right#2}
\newcommand{\st}{\,:\,}
\newcommand{\eps}{\varepsilon}

\newcommand{\I}{[a,b]}

\newcommand{\case}[1] { \left\{ \begin{array}{ll} #1 \end{array} \right. }

\renewcommand{\)}{\right)}
\renewcommand{\[}{\left[}

\newtheorem{theorem}{Theorem}[section]

\newtheorem{proposition}{Proposition}

\theoremstyle{definition}

\newtheorem{remark}{Remark}

\title[Nonsmooth Maximum Principle for mixed and state constraints]
      {A Nonsmooth Maximum Principle for Optimal Control Problems with State and Mixed
      \
       Constraints -- Convex Case}

\author[M. H. A.  Biswas and MdR de Pinho]{}

\subjclass{Primary: 49K15, 49K30; Secondary: 34A60.}
 \keywords{maximum principle, nonsmooth analysis, state constraints.}

 \email{dee08022@fe.up.pt}
 \email{mrpinho@fe.up.pt}

\thanks{The first author is  supported by a   PhD   grant by FCT, Portugal. }

\begin{document}
\maketitle

\centerline{\scshape Md. Haider Ali Biswas }
\medskip
{\footnotesize
 \centerline{Faculdade de Engenharia da Universidade do Porto}
   \centerline{DEEC, Rua Dr. Roberto Frias}
   \centerline{ 4200-465 Porto, Portugal}
} 

\medskip

\centerline{\scshape Maria do Rosario de Pinho}
\medskip
{\footnotesize
 \centerline{Faculdade de Engenharia da Universidade do Porto}
   \centerline{DEEC, Rua Dr. Roberto Frias}
   \centerline{ 4200-465 Porto, Portugal}
}



\begin{abstract}
 Here we derive a nonsmooth maximum principle for
   optimal control problems with both state and mixed
 constraints. Crucial to our development is a convexity assumption on the ``velocity set''.
The  approach consists of applying
 known penalization techniques for state constraints together with  recent results
 for mixed constrained problems.
\end{abstract}

\section{Introduction}

In this paper we develop a nonsmooth maximum principle for optimal control problems with
both pure state and mixed state control constraints in the presence of a convexity assumption. The problem of interest is
 $$(P)\quad
\left\{
 \begin{array}{l}\text{Minimize }  l(x(a),x(b))  \\
\text{subject to}   \\
\begin{array}{l}
\dot{x}(t)  = f(t,x(t),u(t))\quad  \ae ~ t\in [a,b] \\[1mm]
h(t,x(t)) \leq 0 \quad  \hbox{ for all } ~ t\in [a,b] \\[1mm]
      (x(t),u(t))  \in S(t)      \quad    \ae ~ t\in [a,b] \\[1mm]
(x(a),x(b)) \in E.
\end{array}
\end{array}
\right.
$$
The state $x$ and control $u$ are subject to joint, or mixed constraints through
the condition $(x(t), u(t)) \in  S(t)$ where
$t\to S(t)\subset\R^n  \times \R^k$ is a multifunction. The function
$f\colon\R\times\R^{n}\times\R^{k}\rightarrow\R^{n}$
describes the system dynamics and
$h\colon\I\times\R^{n}\rightarrow\R$ is the functional defining the pure state constraint.
Furthermore, the closed set $E\subset \R^{n} \times\R^{n}$
and $l\colon \R^{n}\times\R^{n}\rightarrow\R$
specify the endpoint constraints and cost.

This problem involves  measurable control functions $u$ and absolutely continuous function $x$.
A pair $(x,u)$ is called an admissible  process   if it satisfies  the constraints of the problem with finite  cost.

We say that the process $(\xb,\ub)$ is a \emph{strong local minimum} if, for some $\eps>0$,
it minimizes the cost over admissible processes $(x,u)$ such that
$|x(t)-\xb(t)|\leq \eps\quad \text{for all }t\in [a,b].$
We consider the  \emph{basic hypotheses} on the problem data throughout. They are the following:
$f$ and $L$ are  ${\mathcal L}\times {\mathcal B}^{n+k}$, $S$ is ${\mathcal L}\times {\mathcal B}$, $E$ is closed and $l$ is
locally Lipschitz.

Necessary optimality conditions  for nonsmooth problems with pure state constraints
have been studied systematically for quite some time (see \cite{Vin00} details and references therein).
 On the other hand, problems   with mixed state control constraints,  amply studied
 in a smooth framework (see, for example,
\cite{Arut00}, \cite{Arut05}, \cite{Cla05}, \cite{Dmi93}, \cite{DubMil63},\cite{Hest}, \cite{mil:1},
 \cite{Neu69}, \cite{Zei03})
have received little
 attention. Attempts  to treat  mixed constrained problems with nonsmooth data have been in general timid
 (see, for example, \cite{DL96}, \cite{DLS09}) until quite recently when,
  in  \cite{CdP10}, necessary conditions in the form of a nonsmooth maximum principle were developed.
   However  the literature on nonsmooth maximum principle with both mixed and state constraints has
   been surprisingly  sparse.

In this paper we develop a nonsmmoth maximum principle for problem $(P)$ with both pure  state
and mixed state-control constraints under some convexity assumptions. To achieve our purpose we intertwine  established
approaches  used for state constraints with  up to date developments for problems with mixed constraints.
Indeed, we follow closely the approach developed in \cite{VinPap82} (see also \cite{DFF02} and \cite{DLS09})
 where necessary conditions
for pure state constrained problems  are derived.
Our proofs differ from those in
\cite{DFF02} since we deal not only with state constraints but also with mixed constraints. So applications
of a nonsmooth maximum principle   for mixed constrained problems, derived in
\cite{CdP10} (instead of those in \cite{dPV95}), play a crucial role in our analysis. There is however a price to pay; here
we assume that the solution of $(P)$  is a strong minimum in contrast with  \cite{CdP10}  where a weaker notion of minimum,
that of
 \emph{local minimum of radius R}  is used (in this respect see also \cite{memoirs}).
Also we need to strengthen the  hypotheses in comparison with those in \cite{CdP10}.
The convexity assumption we impose on this paper may be seen as
 a major hindrance to some applications.
 Although  this  assumption can be successfully
removed following the lines of \cite{DFF05} we opt, for the sake of simplicity, to  report that work
elsewhere together with a discussion of the hypotheses
and illustration of applications.
%
%
%

\section{Preliminaries}

For $g$ in $\R^m$, inequalities like $g\leq 0$
are interpreted componentwise. 
Here and throughout, $\B$ represents the closed unit ball centered
at the origin regardless of the dimension of the underlying space
and  $\abs{\,\cdot\,}$ the Euclidean norm or the induced matrix
norm on $\R^{p\times q}$.
  The {\em Euclidean
distance function} with respect to a given set $A \subset \R^{m}$
is
$$
d_{A}\colon\R^{k}\rightarrow\R, \qquad
y \mapsto d_{A}(y) = \mbox{inf }\Set{ \abs{y-x} \st x \in A}.
$$

A function $h\colon\I\to\R^p$ lies in $W^{1,1}(\I;\R^p)$
if and only if it is absolutely continuous;
in $L^1(\I;\R^p)$ iff it is integrable; and
in $L^{\infty}(\I;\R^p)$ iff it is essentially bounded. The norm  of $L^\infty(\I;\R^p)$ is $\norm{\cdot}_\infty$.

\smallskip

The space $C^*([a,b];\mathbb{R}) $ is  the topological dual   of the space of continuous functions $C([a,b];\mathbb{R}) $.
 Elements of  $C^*([a,b];\mathbb{R}) $ can be identified with finite regular measures on the Borel subsets of [a,b].
 The set of elements in  $C^*([a,b];\mathbb{R}) $ taking nonnegative values on nonnegative-valued functions in
 $C([a,b];\mathbb{R}) $ is denoted by
 $C^\oplus([a,b];\mathbb{R})$. The norm in $C^\oplus([a,b];\mathbb{R})$, $|\mu|$  , coincides with the total variation of
$\mu, ~\int_{[a,b]}\mu(ds)$. The support of a measure $\mu$, written  $\hbox{ supp} \{\mu\} $,
is the smallest closed set  $A\subset \I$ such that for any relatively open subset $B\subset [a,b]\backslash A $  we have
$\mu(B)=0$ .
\smallskip



We make use of standard concepts from nonsmooth analysis.
Let $A\subset \R^k$ be a closed set with $\xb\in A$.
The  {\em proximal normal cone to $A$ at $\xb$} is denoted by
$N_{A}^P (\xb)$, while $N_{A}^L (\xb)$ denotes the {\em limiting normal cone }
and $N_{A}^C (\xb)$ is the {\em Clarke normal cone}.

Given a lower semicontinuous function
$f\colon\R^k\rightarrow\R\cup\Set{+\infty}$
and a point $\xb\in \R^{k}$ where $f(\xb) < +\infty $,
$\partial^L f(\xb)$ denotes the
{\em limiting subdifferential} of $f$ at $\xb$.
When the function $f$ is Lipschitz continuous near $x$,
the convex hull of the limiting subdifferential, $\mbox{co }\partial^L f(x)$,
coincides with the {\em (Clarke) subdifferential} $\partial^C f(\xb)$.
For details on such nonsmooth analysis concepts, see for example \cite{Cla83,clsw98,Mor05,Vin00}.

\section{Auxiliary Results}
\label{FC}

In this section we present a simplified version of one of the main results in \cite{CdP10} that
 will be of importance in the forthcoming developments.


Take a fixed interval $[a,b]$ and a set $S$ of $[a,b]\times \R^n\times \R^k$.
Define
\begin{equation}\label{setS}S(t):=\left\{(x,u): (t,x,u)\in S\right\} \text{ for all } t\in\I.\end{equation}
Assume  for the time being that  $E\subset \R^n\times\R^n$ and $l:\R^n\times \R^n \to \R$. Consider the following problem:
 $$(C)~
\left\{
 \begin{array}{l}\text{Minimize }  l(x(a),x(b)) +\disp\int_a^b L(t,x(t), u(t))dt \\
\text{subject to}   \\
\begin{array}{rcl}
\dot{x}(t)  &=   & f(t,x(t),u(t)) ~ \ae ~ t\in [a,b] \\
      (x(t),u(t))  &\in & S(t)           ~ \ae ~ t\in [a,b] \\
(x(a),x(b)) &\in & E.
\end{array}
\end{array}
\right.
$$
where $L:[a,b]\times \R^n\times\R^k\to \R$.

For some $\eps>0$\footnote{The $\eps$ here can be taken to be equal to the parameter defining the strong local minimum.}
define
$$S^\epsilon_*(t)=\left\{(x,u) \in S(t):~ |x-\xb(t)|\leq  \epsilon\right\}.$$

In generic terms we assume that a function $\phi(t,x,u)$ satisfies $[L_*^\epsilon]$ if:

\noindent $\mathbf{[L_*^\epsilon]}$ There exist constants $k_x^\phi$ and $k_u^\phi$ such that for almost every $t\in [a,b]$ and every
$(x_i,u_i) \in S^\epsilon_*(t)$ ($i=1,2$) we have
$$|\phi(t,x_1,u_1)-\phi(t,x_2,u_2)|\leq k_x^\phi|x_1-x_2|+k_u^\phi |u_1-u_2|.$$

If this assumption is imposed on $f$, then the Lipschitz constants are denoted by $k_x^f$ and $k_u^f$.
As for $S(t)$ we consider the following \textbf{bounded slope condition}:

\noindent  $\mathbf{[BS_*^\epsilon]}$ There exists a  constant $k_S$ such that for almost every $t\in [a,b]$ the following condition holds
    $$(x,u) \in S^\epsilon_*(t),~(\alpha,\beta) \in N_{S(t)}^P(x,u) \Longrightarrow |\alpha|\leq k_S|\beta|.$$

The two previous hypotheses are strengthening of the analogous hypotheses in \cite{CdP10}.
For the sake of uniformity and the analysis in the forthcoming sections we need to position an extra
hypothesis on the set  $S^\epsilon_*(t)$. We assume that:

\noindent  $\mathbf{[CS_*^\epsilon]}$ The set $S^\epsilon_*(t)$  is closed and there exists an integrable function $c$ such that
for almost every $t\in [a,b]$  the following holds
$$S^\epsilon_*(t) \hbox{ is closed and }(x,u) \in S^\epsilon_*(t)~\Longrightarrow |(x,u)|\leq c(t).$$
\smallskip

We observe that although $[CS_*^\epsilon]$ is a strong assumption it is nevertheless of  importance in our future development. Necessary conditions of optimality for $(C)$ are given by the following theorem:

\begin{theorem}\label{FC1} ({\rm adaption of   Theorem 7.1 in \cite{CdP10}})
Let $(\xb,\ub)$ be a strong local minimum for problem~$(C)$. Assume that
the basic hypotheses, that $f$ and $L$  satisfy $[L_*^\epsilon]$ and that  $[BS_*^\epsilon]$  and $[CS_*^\epsilon]$
hold.

 Then there exist an absolutely continuous function
$p\colon\I\to\R^n$,
  and a scalar
$\lambda_0 \ge 0$ such that
\beeq \label{3:t11}
 & (p(t),  \lambda_0)\neq 0 \quad \forall~ t\in [a,b],\qquad\qquad \\[2mm]
\label{3:t12}
&  (-\dot{p}(t),0)\in  \partial_{x,u}^C \left[\langle p(t),\bar{f}(t)\rangle -\lambda_0 \bar{L}(t)-K|p(t)| \bar{d}_{S(t)}(t)\right]
~ \text{a.e. }
\\[2mm]
\label{3:t13}
&  (\xb(t),u)\in S(t)\Longrightarrow  \langle p(t),f(t,\xb(t),u)\rangle -\lambda_0L(t,\xb(t),u)\\
\nonumber
& \leq \langle p(t),f(t,\xb(t),\ub(t))-\lambda_0 L(t,\xb(t),\ub(t))\rangle
~\text{a.e. }
\\[2mm]
\label{3:t14}
& (p(a), -p(b)) \in  N_{E}^L(\xb(a), \xb(b)) + \lambda_0\d l(\xb(a),
\xb(b)),
\eeeq
where $\bar{f}(t)$ and $\bar{L}(t)$ represent the function evaluated at $(t,\xb(t),\ub(t))$, $\bar{d}_{S(t)}(t)$
is the distance function to $S(t)$ evaluated at $(\xb(t),\ub(t))$ and $K$ in (\ref{3:t12}) is
a constant  depending only on
 $k_x^f,~k_x^L, ~k_u^f,~k_u^L$ and  $k_S$.
\end{theorem}

\section{The Convex Case} We now turn to problem $(P)$. For this problem we derive
a nonsmooth maximum principle under the following convexity assumption  on the ``velocity set'':

\begin{itemize}
\item[\bf{[C]}] The velocity set
$\left\{ v\in \R^n:~v=f(t,x,u),~u\in S(t,x)\right\}$
is convex for all $t\in [a,b]$.
\end{itemize}

Furthermore we need to  impose two more hypotheses on the data of our problem,
one related to the state constraint and another to  mixed constraints.

\begin{itemize}
\item[\bf{[H1]}] For all $ x \in \xb(t) + \eps\B$ the function $t \to h(t, x)$ is continuous and
there exists a scalar $k_h>0$  such that the function $x \to h(t, x)$ is
Lipschitz of rank $k_h$ for all $t \in \I$.
\item[\bf{[H2]}] For almost every $t\in [a,b]$ the following condition holds:
for all $u\in S(t,\xb(t))$ and all sequence $x_n\to \xb(t)$   there exists a sequence $u_n \in S(t,x_n)$ such that $u_n\to u$.
\end{itemize}
In the above  the set $S(t,x)$ is defined as
$$S(t,x)=\left\{u:~(x,u)\in S(t)\right\}$$ where $S(t)$ is as in (\ref{setS}).
For a discussion on the need to impose continuity of $t\to h$ see \cite{DFF02}.
Hypothesis [H2] asserts the lower semi-continuity of the multifunction $x\to S(t,x)$ (for definition and properties see
\cite{Frank}).

Assume the basic assumptions. Also suppose that $f$ satisfies $[L_*^\epsilon]$ and that both
  $[BS_*^\epsilon]$ and $[CS_*^\epsilon]$ hold. Under these assumptions we  note
  for future use that the following conditions are
satisfied:
\begin{equation}\label{Lip}
|f(t,\xb(t),u)-f(t,\xb(t),\ub(t))|\leq   k_u^f|u-\ub(t)|\text{ for all }u \in S(t,\xb(t))\text{ a.e. }t
\end{equation}
 for all $u\in S(t,\xb(t))$ a.e. $t\in [a,b]$ and there exists an integrable
function $k$ such that
\begin{equation}\label{boundoff} | f(t,\xb(t),u)|\le k(t) \text{ for all }u \in S(t,\xb(t))\text{ a.e. }t .
\end{equation}
%
%

Before proceeding we need to define the following subdifferential
\begin{equation}\label{barpartialh}
  \bar{\partial}_x h(t,x) :=\co \{\lim \xi_i\ :\ \xi_i \in \partial_x h(t_i,x_i), (t_i,x_i) \to (t,x) \}.
\end{equation}

We are now in position to state our main result.

\begin{theorem}\label{DFF}
Let $(\xb,\ub)$ be a strong local minimum for problem~$(P)$. Assume that
the basic hypotheses,  [C], [H1], [H2], $[BS_*^\epsilon]$  and $[CS_*^\epsilon]$
hold and that $f$ satisfies $[L_*^\epsilon]$.
 Then there exist an absolutely continuous function
$p\colon\I\to\R^n$, an integrable function $\gamma:\I \to \R^n
   $,
     a  measure $\mu \in C^\oplus([a,b];\mathbb{R})$,
  and a scalar
$\lambda_0 \ge 0$ such that
$$\begin{array}{lc}
\hspace{-0.5cm}\rm{(i)}  & \mu \{[a,b]\} + ||p||_{\infty}+\lambda_0 > 0,\qquad\qquad \\[2mm]
\hspace{-0.5cm}\rm{(ii)}
 & (-\dot{p}(t),0)\in \partial^C_{x,u}\langle q(t), f(t,\xb(t),\ub(t))\rangle- N_{S(t)}^C(\xb(t),\ub(t))~~ \text{a.e.},
\\[2mm]
\hspace{-0.5cm}\rm{(iii)}
& (\xb(t),u)\in S(t)\Longrightarrow
 \langle q(t),f(t,\xb(t),u)\rangle \leq \langle q(t),f(t,\xb(t),\ub(t))\rangle
~~\text{a.e.},
\\[2mm]
\hspace{-0.5cm}\rm{(iv)}
& (p(a), -q(b)) \in  N_{E}^L(\xb(a),\xb(b)) + \lambda_0\d l(\xb(a),
\xb(b)),\\[2mm]
\hspace{-0.5cm}\rm{(v)}   & \gamma(t) \in \bar{\partial} h(t,\xb(t)) \quad \mu \mbox{-} \text{a.e.},\\[2mm]
\hspace{-0.5cm}\rm{(vi)}   & \mathrm{supp}\{ \mu \} \subset
     \left\{ t \in [a,b] : h\(t,\xb (t) \)=0  \right\},
\end{array}$$
where
 \begin{equation}
    \label{qqq}
    \begin{split}
    q(t)=\case{ p(t)+\int_{[a,t)} \gamma( s) \mu (ds) &  t \in [a,b)\\
                p(t)+\int_{[a,b]} \gamma( s) \mu (ds) &  t =b.}
 \end{split}
  \end{equation}
\end{theorem}

\begin{remark}The proof shows  that we prove a sharper
form of (ii):
$$(-\dot{p}(t),0)\in \partial^C_{x,u}\langle q(t), f(t,\xb(t),\ub(t))\rangle- K|q(t)|\partial_{x,u}^C d_{S(t)}(\xb(t),\ub(t)).$$

\end{remark}

\begin{remark}
It is also easy to deduce from the proofs that when assumption [H2] is not imposed, a ``weaker'' version of
the necessary conditions for $(P)$  (in the vein  of \cite{dPV95}) can be obtained:
 all the conclusions but
(iii) (the Weierstrass condition) hold.
\end{remark}

We derive Theorem~\ref{DFF} in two main stages.
In the first stage we establish the validity of the theorem to the following problem

$$
\mathrm{(Q)}\hspace{.2in}\left\{
 \begin{array}{l}\text{Minimize }  l(x(b))  \\
\text{subject to}   \\
\begin{array}{rclrr}
\dot{x}(t)  &=   & f(t,x(t),u(t)) & \quad\ \ \quad\ae & t\in \I \\
         (x(t),u(t))  &\in & S(t) & \quad\ \ \quad\ae & t\in \I \\
     h(t,x(t))  &\leq & 0        & \quad\ \ \text{for all }& t\in \I \\
(x(a),x(b))&\in & \{x_a\}\times E_b.             &                   &
\end{array}
\end{array}
\right.
$$
Problem $(Q)$ is a special case of $(P)$ in which
$E=\{x_a\}\times E_b$ and $l(x_a,x_b)=l(x_b)$.

\begin{proposition}
\label{FEC}
Let $(\xb,\ub)$ be a strong local minimum for problem~$(Q)$. Assume
the basic hypotheses, [C], [H1],  $[BS_*^\epsilon]$  and $[CS_*^\epsilon]$ hold, that $E_b$ is closed
 and that $f$ satisfies $[L_*^\epsilon]$.
 Then there exist an absolutely continuous function
$p\colon\I\to\R^n$, an integrable function $\gamma:[a,b] \to \R^n
   $,
      a  measure $\mu \in C^\oplus([a,b];\mathbb{R})$,
  and a scalar
$\lambda_0 \ge 0$ such that
\beeq \label{p:t11}
 & \mu \{[a,b]\} + ||p||_{\infty}+\lambda_0 > 0,\qquad\qquad \\[2mm]
\label{p:t12}
 & (-\dot{p}(t),0)\in \partial^C_{x,u}\langle q(t), f(t,\xb(t),\ub(t))\rangle- N_{S(t)}^C(\xb(t),\ub(t))~~ \text{a.e. }
\\[2mm]
\label{p:t13}
& (\xb(t),u)\in S(t)\Longrightarrow
 \langle q(t),f(t,\xb(t),u)\rangle \leq \langle q(t),f(t,\xb(t),\ub(t))\rangle
~~\text{a.e. }
\\[2mm]
\label{p:t14}
& -q(b) \in  N_{E_b}^L(\xb(b)) + \lambda_0\d l(
\xb(b)),\\[2mm]
\label{p:t15}  & \gamma(t) \in \bar{\partial} h(t,\xb(t)) \quad \mu \mbox{-} \ae,\\[2mm]
\label{p:t16}& \mathrm{supp}\{ \mu \} \subset
     \left\{ t \in [a,b] : h\(t,\xb (t) \)=0  \right\},
\eeeq
where $q$ is as in (\ref{qqq}).

 \end{proposition}

\section{Proof of Proposition \ref{FEC}}
We now proceed  proving Proposition \ref{FEC}.
Observe that $u(t)\in S(t,x(t))$ is equivalent to $(x(t),u(t)) \in S(t)$.

The local minimality of $(\xb,\ub)$  provides some $\varepsilon>0$.
By reducing this constant if necessary, we can also rely on the hypotheses.
The proof breaks into several steps.

\noindent {\bf  Step 1:} \emph{Penalize  state-constraint violation.}

Define the following problem for each $i\in \N$:
$$(Q_i)
\quad
\left\{
 \begin{array}{l}\text{Minimize }  l(x(b)) + i \displaystyle \int_a^b h^+(t,x(t))~dt \\
\text{subject to}   \\
\begin{array}{l}
 \dot{x}(t)  = f(t,x(t),u(t))\quad  \ae ~ t\in [a,b] \\[1mm]
      (x(t),u(t))  \in S(t)      \quad    \ae ~ t\in [a,b] \\[1mm]
(x(a),x(b)) \in \{x_a\}\times E_b,
\end{array}
\end{array}
\right.
$$
where $
 h^+(t,x):=\max \{ 0, h(t,x)\}.$ This differs from~$(Q)$  by shifting the state constraint
into the objective function.

Following the approach in \cite{VinPap82}  (see also \cite{DFF02})
let us temporarily assume that penalization is effective, i.e.,
\begin{itemize}
\item[\bf{[IH]}] $\disp\lim_{i\to \infty}
\inf\{P_i\} = \inf\{P\}.$
\end{itemize}
We will justify this assumption later.

\noindent {\bf  Step 2:} \emph{Application of Ekeland's theorem.}

Let $W$ denote the set of measurable functions
$u\colon\I\to\R^k$   for which there  exists an
absolutely continuous function $x$ such that
$
\dot{x}(t)=f(t,x(t),u(t))$, $(x(t),u(t)) \in  S(t)$,    for almost every $t\in \I$,
   $ x(t) \in   \xb(t)+\varepsilon \B$ for all $t\in \I$, $x(a)=x_a$ and $x(b)\in E_b$.
We provide $W$ with the metric
$
 \Delta(u,v):=\parallel u-v\parallel_{L_1}
$
and define $J_i\colon W\to\R$ using the arc $x$ mentioned above:
$$
J_i(u):=l(x(b)) + i\disp\int_a^b h^+(t,x(t))\,dt.
$$
It is a simple matter to check that $(W,\Delta)$
is a complete metric space in which the functional
$J_i\colon W\rightarrow\R$ is continuous (see \cite{Cla83}).
Moreover, problem~$(Q_i)$ above is closely related to
the abstract problem
$$
(R_i)\quad\left\{\begin{array}{ll}
\text{Minimize }  & J_i(u)        \\
\text{subject to} & u \in W.
\end{array}\right.
$$
Clearly $(\ub,\xb(b))$ is admissible for $(R_i)$, with
$
J_i(\ub)
= l(\xb(b))
= \mathrm{inf} ~P
$
since for all  $t\in [a,b]$, $h^+(t,\xb(t))= 0$ .
Let
$
\eps_i
= J_i(\ub)
  - \mathrm{inf}~P_i.
$
We have $\eps_i\geq 0$ and, taking into account [IH], $\eps_i\to 0$.
Ekeland's variational principle (see \cite{Vin00}) applies. It
asserts the existence of $u_i\in W$ such that
\beeq\label{norma11} \parallel u_i-\ub\parallel_{L_1}
\leq \sqrt{\eps_i} \eeeq and $u_i$ minimizes over $W$ the
perturbed cost functional \beeq\label{perturbed11} u \mapsto
J_i(u)+\sqrt{\eps_i}\parallel u_i-\ub\parallel_{L_1}. \eeeq
Let $x_i$ be the trajectory corresponding to $u_i$.

\noindent{\bf Step 3:}
{\em Study optimality conditions for the perturbed problem.}

In control-theoretic notation, our work with
Ekeland's Theorem shows
that the process $(x_i,u_i)$ solves the following optimal control problem:
$$
(D_i)\quad\left\{
\begin{array}{l}
\text{Minimize }  l(x(b)) + i\disp\int_a^b h^+(t,x(t))\,dt+ \sqrt{\eps_i}\int_a^b \abs{u(t)-u_i(t)}~ dt \\
\text{subject to}   \\
\begin{array}{rclrr}
\dot{x}(t) & =   & f(t,x(t),u(t))    & \qquad\ae &t\in \I \\
      (x(t),u(t)) & \in & S(t)       & \qquad\ae &t\in \I \\
      x(t) & \in & \xb(t)+\varepsilon \B & \qquad\text{ for all} &t\in \I \\
x(a)& = & x_a    & &        \\
    x(b) & \in  & E_b.              & &
\end{array}
\end{array}
\right.
$$
Since $\eps_i\to 0$ (by [IH]) it follows from (\ref{norma11}) that
 $u_i\to \ub$ strongly. We can then arrange by subsequence extraction, if necessary, that
 $u_i\to \ub$ almost everywhere.   We can further deduce that $x_i\to \xb$ uniformly. By
discarding initial terms of the sequence we can guarantee that
$(x_i,u_i)$ is a local minimum
for a variant of problem $(D_i)$ obtained by dropping the constraints
 $ x(t)  \in \xb(t)+\eps\B$.
 We now fix our attention in the related subsequence of problems without relabeling.

\smallskip

Theorem \ref{FC1} applies to $(D_i)$. It provides an absolutely continuous function
$p_i$ and a scalar $\lambda_i\geq 0$ such that
\begin{equation}\label{nc0di}
(p_i(t),\lambda_i)\neq  0 \text{ for all }t,
\end{equation}
\begin{equation}\label{nc1di}
\begin{split}
 (-\dot p_i(t),0) \in ~& \partial_{x,u}^C \left\{ \langle p_i(t),f(t,x_i(t),u_i(t))\rangle -i\lambda_i h^+(t,x_i(t))\right.\\
& \left. -\sqrt{\eps_i}\lambda_i \abs{u(t)-u_i(t)} -K|p_i(t)|d_{S(t)}(x_i(t),u_i(t))\right\}
~\text{a.e.}
\end{split}
\end{equation}
\begin{equation}\label{nc2di}
\begin{split}
 (x_i(t),u)\in S(t)&\Longrightarrow\\
 \langle p_i(t),f(t,x_i(t),u)\rangle -\sqrt{\eps_i}\lambda_i\abs{u(t)-u_i(t)}& \leq \langle p_i(t),f(t,x_i(t),u_i(t))\rangle
~\text{a.e.}
\end{split}
\end{equation}
\begin{equation}\label{nc3di} -p_i(b) \in N^L_{E_b}(x_i(b))
 +\lambda_i\partial^L l(x_i(b))\end{equation}

These conditions have consequences we now seek to express in terms
of the original problem~$(Q)$. 

Apply Clarke's sum rule~\cite{Cla83}
to~(\ref{nc1di}) and take into accounts the properties of the subdifferentials of the distance function. We deduce that there exist measurable functions $\xi_i,~\zeta_i,~\gamma_i,~e_i,~\phi~_i$ and $\varphi_i$ such that for almost every $t$ in $\I$,
\beeq
\label{pi} & (\xi_i(t),\zeta_i(t)) \in \partial^C_{x,u} f(t,x_i(t),u_i(t)),\\
\label{hi} & (\gamma_i(t),0) \in \partial^C_{x,u} h^+(t,x_i(t)),\\
\label{ui}  & e_i(t) \in \R^k \text{ such that } |e_i|\leq 1,\\
\label{di} & (\phi_i(t),\varphi_i(t)) \in \partial^C_{x,u}  d_{S(t)}(x_i(t),u_i(t)),\quad |(\phi_i(t),\varphi_i(t))|\leq 1
\eeeq
such that

\begin{equation}
\label{piequal}
-\dot p_i(t) = p_i(t)\xi_i(t)-i\lambda_i\gamma_i(t)-K|p_i(t)|\phi_i(t),
\end{equation}
\begin{equation}
\label{zeroequal}
0 = p_i(t)\zeta_i(t)-\sqrt{\eps_i}\lambda_i e_i-K|p_i(t)|\varphi_i(t).
\end{equation}

To simplify this further,
let $h_0(t,x)=0$ and $h_1(t,x)=h(t,x)$ so that
$$h^+(t,x)=\max\Set{h_j(t,x)\st j=0,1}.$$
Then for each fixed $t$, Clarke's Max Rule~\cite{Cla83} says
$$\begin{array}{rl}
\partial_{x,u}^C  h^+(t,x_i(t)) &
\\
&\disp\hskip -6em
\subseteq \text{co}\disp\cup_{j=0}^1\Set{\partial_{x,u}^C h_j(t,x_i(t))
                         \st h_j(t,x_i(t)) = h^+(t,x_i(t))}.
\end{array}$$
Clearly $\partial_{x,u}^C h_0\equiv\Set{(0,0)}$,
so a typical element of the right side has the form
$
 \alpha_i\lr(){\gamma_i,~0},
$
where $(\gamma_i,0)\in\partial_{x,u}^C h_1(t,x_i(t))$
and $\alpha_i$ is chosen from
$$
\Sigma_i(t)
\!= \!\Set{\!\alpha~\in [0,1],~\alpha = 0~\text{if}\ h_1(t,x_i(t))<h^+(t,x_i(t))\!}.
$$
Tracking
these dependencies leads to the following expansion of (\ref{piequal}) and (\ref{zeroequal}):
\begin{equation}\label{13}
\begin{split}
-\dot
p_i(t) = &~ p_i(t)\xi_i(t)-i\lambda_i\alpha_i(t)\gamma_i(t)-K|p_i(t)|\phi_i(t), \\[2mm]
0= & ~p_i(t)\zeta_i(t)-\sqrt{\eps_i}\lambda_i e_i(t)-K|p_i(t)|\varphi_i(t).
\end{split}
\end{equation}

We now introduce the measure $\mu_i\in C^*(\I;\R)$:
$$\disp\int_B d\mu_i(t)=\disp\int_B i\lambda_i \alpha_i(t)dt$$
for every Borel set $B\subset \I$.
Define  $\pi_i\in C^*(\I;\R)$ as
$d\pi_i(t)=\dot p_i(t)dt.$
Then, from (\ref{13}) we get
\begin{equation}\label{ajointequation}
\begin{split}-\disp\int_B d\pi_i(t)=~& \disp\int_B\big( p_i(t)\xi_i(t)-K|p_i(t)|\phi_i(t)\big)dt-\disp\int_B
\gamma_i(t)d\mu_i(t)\\
0~=~& \disp\int_B \big(p_i(t)\zeta_i(t)-\sqrt{\eps_i}\lambda_i e_i(t)-K|p_i(t)|\varphi_i(t)\big)dt
\end{split}
\end{equation}
and we have \begin{equation}
\label{pipi}p_i(t)=b_i+\disp\int_{[a,t)}d\pi_i(t)~\text{ for all } t\in (a,b],
\end{equation}
for every Borel set $B$. Here $b_i=p_i(a)$.
Taking  (\ref{nc3di}) into account we have
\begin{equation}\label{transversality}-b_i-\disp\int_{[a,b]}d\pi_i(t)\in N^L_{E_b}(x_i(b))
 +\lambda_i\partial^L l(x_i(b)).\end{equation}
Since $\alpha_i(t) \in \Sigma_i(t)$, we have $\mu_i \in C^\oplus(\I;\R)$ and this measure has support in
$\left\{t\in \I: h(t,x_i(t))\geq 0\right\}$.

Since, by (\ref{nc0di}) $b_i$ and $\lambda_i$ are not both zero, we may conclude, after rescaling, that
\begin{equation}
\label{nontriviality}
|b_i|+|\mu_i|+\lambda_i=1.
\end{equation}

\noindent{\bf Step 4:}
{\em Take limits.}
Our first steps has dealt with fixed $i\in \N$. We now consider the case when $i\to \infty$.
Recall that the sequence $x_i$ converges uniformly to $\xb$ and $u_i\to \ub$ almost everywhere.

Under the hypotheses and appealing to Gronwall's inequality we deduce the existence of a constant $K_1$ such that
$|\pi_i|\leq K_1.$
It follows from (\ref{pipi}) and (\ref{nontriviality}) that
$|p_i(t)|\leq K_1+1.$
We now deduce from the above that
$$\pi_i\to \pi ~~\text{weakly}^*$$
for some measure $\pi$. Consequently
$$|\pi_i|\to |\pi|.$$
Turning again to  (\ref{nontriviality}) we may arrange that
$b_i\to b,\quad \lambda_i\to \lambda,\quad \mu_i\to \mu~~ \text{weakly}^*$
for some $b\in\R^n$, $\lambda\geq 0$ and some measure $\mu$. Also we have
$|\mu_i|\to |\mu|$ and
$$|b|+|\lambda|+|\mu|=1.$$
With the above and appealing to Lemma 4.3 in \cite{VinPap82} we can now conclude that
there exists some subsequence such that
$p_i(t)\to q(t)~~\text{ a.e. }$
where $q$ is now a function of bounded variation defined as
$$q(t):=b+\disp\int_{[a,t)} d\pi\quad \text{and }\quad b_i+\disp\int_{[a,t)} d\pi_i\to b+\disp\int_{[a,t)} d\pi.$$

Under the hypotheses we deduce from (\ref{pi}) that
$|(\xi_i(t),\zeta_i(t))|\leq \max\{k_x^f,k_u^f\}$ a.e.
Dunford-Pettis Theorem (see for example \cite[Theorem 2.51]{Vin00})
asserts existence of a subsequence converging weakly in the $L^1$ topology to some function $(\xi,\zeta)$ such that
$\xi,~\zeta \in L^1$. Taking into account
(\ref{ui}) and (\ref{di}) we deduce in the same way that
$ e_i \to e,\quad
(\phi_i,\varphi_i)\to (\phi,\varphi)$ for some
$e,~\phi,~\varphi \in L^1$ where the convergent is understood in the weak $L^1$ topology. Upper semi-continuity properties of the subdifferentials asserts that
(\ref{pi})--(\ref{di}) hold when we remove the indexes $i$.

Observe that  $\partial_{x}^C h(t, x) \subset \bar{\partial}_x h(t,x)$ (see (\ref{barpartialh}) for definition of $\bar{\partial}_x h(t,x)$) and that  $\bar{\partial}_x h(t,x)$ is of closed graph for any $i$. It follows from \cite[Lemma 4.3]{VinPap82}
that there exists   a Borel measurable, $\mu$-integrable function $\gamma$
such that
$\gamma(t) \in \bar{\partial}_x h(t,\xb(t))\quad \mu-\text{a.e.}$
This is (\ref{p:t15}) of the proposition.

We now turn to (\ref{transversality}). The properties of limiting normal cones and limiting subdifferential assert
that
\begin{equation}\label{transversality2}
-b-\disp\int_{[a,b]}d\pi(t)\in N^L_{E_b}(\xb(b))
 +\lambda\partial^L l(\xb(b)).
 \end{equation}

We  concentrate on the support of the measure $\mu$. Mimicking the arguments in \cite{DFF02}
it is a simple matter to see that
$\text{supp}\{\mu\}\subset\left\{t\in[a,b]:~h(t,\xb(t))=0\right\}.$
This is  conclusion (\ref{p:t15}) of the proposition.

By [H2] we deduce from (\ref{nc2di}) that
 \begin{equation}\label{weierstrass}
 \langle q(t),f(t,\xb(t),u)\rangle
 \leq \langle q(t),f(t,\xb(t),\ub(t))\rangle.
 \end{equation}

Next we focus on (\ref{ajointequation}). Lemma 4.3 in \cite{VinPap82} and our conclusions above
assert that
$$\begin{array}{rcl}-q(t)+b& = & \disp\int_{[a,t)} \Big(q(s)\xi(s)-K|q(s)|\phi(s)\Big)ds-\disp\int_{[a,t)}\gamma(s)d\mu(s)\\[3mm]
0 & = & \disp\int_{[a,t)} \Big(q(s)\zeta(s)-K|q(s)|\varphi(s)\Big)ds
\end{array}$$
Define now the function
$p(t):=q(t)-\disp\int_{[a,t)} \gamma(s)d\mu(s).$
From the above we now obtain the conclusions of the proposition. Observe that
(\ref{p:t11}) follows from $|b|+|\lambda|+|\mu|=1.$

\noindent{\bf Step 5:}
{\em Show that [C] implies [IH].}
We omit the details since the conclusion can be obtained adapting the arguments in \cite{DFF02}.

\section{Sketch of the Proof of Theorem \ref{DFF}}
The proof comprises three stages.
We omit  the details.
We first extend Proposition \ref{FEC} to problems where
$x(a) \in E_a$
and $E_a$ is a closed set.
This is done following the lines  in the end of the proof of Theorem 3.1 in \cite{VinPap82}.
Thus we obtain necessary conditions when $(x(a),x(b))\in E_a\times E_b$.
Next we consider the case when the cost is $l=l(x(a),x(b))$. This is done using the technique
in Step 2 of section 6 in \cite{DLS09}.
And finally,  following the approach in  section 6 in \cite{DLS09}, we
derive necessary conditions when
$(x(a),x(b)) \in E$ and $E$ is a closed set.
This completes the proof.

\section*{Acknowledgments}
Md. Haider Ali Biswas would like to thank the Doctoral Program PDEEC,
 Faculdade de Engenharia da Universidade do Porto  for
the grant he received
when starting  this work.
The financial support of FCT Project
PTDC/EEA-CRO/116014/2009,
is also acknowledged.


\medskip
Received July 2010; revised April 2011.
\medskip

\end{document}